\documentclass[leqno,11pt]{article}
\usepackage[spanish]{babel}
\usepackage[ansinew]{inputenc}
\usepackage{amscd}
\usepackage{amsthm}
\usepackage{amssymb}
\usepackage{amsmath}
\usepackage{graphics}
\usepackage{graphicx}
\usepackage{verbatim}
\usepackage{color}
\newtheorem{mydef}{Definition}
\newtheorem{myexa}{Example}

\newtheorem{mytheo}{Theorem}
\newtheorem{mylemma}{Lemma}

\newtheorem*{myprem}{{\em Remarks}}

\newcommand{\cind}{\perp\hspace{-1.25ex}\perp}

\newcommand{\pt}{\mbox{$\succ$\hspace{-1ex}$\longrightarrow$}}

\decimalpoint

\begin{document}

\begin{center}
	{\sc A Note on Conditional Expectation for Markov Kernels}\vspace{2ex}\\
	A.G. Nogales\vspace{2ex}\\
	Dpto. de Matem\'aticas, Universidad de Extremadura\\
	Avda. de Elvas, s/n, 06006--Badajoz, SPAIN.\\
	e-mail: nogales@unex.es
\end{center}
\vspace{.4cm}
\begin{quote}
	\hspace{\parindent} {\small {\sc Abstract.} A known property of conditional expectation is extended to the framework of Markov kernels. Its meaning in terms of densities is provided. Some examples located in the field of clinical diagnosis are presented to delimit the main result of the paper. }
\end{quote}

\vfill
\begin{itemize}
	\item[] \hspace*{-1cm} {\em AMS Subject Class.} (2010): {\em Primary\/} 60Exx
	{\em Secondary\/} 60J35
	\item[] \hspace*{-1cm} {\em Key words and phrases:} conditional expectation, Markov kernel.
\end{itemize}

\newpage

\section{Introduction and basic definitions}

A known property of conditional expectation states that, given an integrable real random variable $X$ and two sub-$\sigma$-fields $\mathcal B_i$, $i=1,2$, 
$$E(X|\mathcal B_1\vee\mathcal B_2)=E(X|\mathcal B_1),
$$
provided that $\mathcal B_2$ is independent of the $\sigma$-field $\mathcal B_1\vee\sigma(X)$ generated by $\mathcal B_1\cup\,\sigma(X)$, where $\sigma(X)$ is the $\sigma$-field generated by $X$. 
See Williams (1991, p. 88, 9.7.(k)) for instance.

In terms of random variables the result reads as follows: if $Y,Z$ are random variables such that $Z$ is independent of $(X,Y)$, then 
$$E(X|Y,Z)=E(X|Y).
$$

It is the main aim of this note to obtain a generalization of this result  for Markov kernels. Some examples, set within the framework of clinical diagnosis, are presented to delimit our main result.

The concepts presented in this section can be found in Heyer (1982) (see also Dellacherie and Meyer (1988)) or in previous paper by the author, and therefore they will be exposed very briefly, even at risk of being somewhat dense.  
 
However the usual notations in this area have been modified. 
It is well known that the concept of  Markov kernel is an extension of the concept of random variable (and also of the concept of $\sigma$-field) and the notation to be used for operations with Markov kernels,  the same  as for random variables, tries to highlight this analogy.

In the next,  $(\Omega,\mathcal A)$, $(\Omega_1,\mathcal A_1)$,
and so on, will denote measurable spaces. A random variable is a map
$X:(\Omega,\mathcal A)\rightarrow (\Omega_1,\mathcal A_1)$ such that
$X^{-1}(A_1)\in\mathcal A$, for all $A_1\in\mathcal A_1$. Its
probability distribution (or, simply, distribution) $P^X$ with
respect to a probability measure $P$ on $\mathcal A$ is the image
measure of $P$ by $X$, i.e., the probability measure on $\mathcal
A_1$ defined by $P^X(A_1):=P(X^{-1}(A_1))$. We will write $\times$
instead of $\otimes$ for the product of $\sigma$-fields or measures. $\mathcal R^k$ will denote the Borel $\sigma$-field on $\mathbb R^k$.

\begin{mydef}\rm  (i) (Markov kernel) A Markov kernel
	$M_1:(\Omega,\mathcal A)\pt   (\Omega_1,\mathcal A_1)$ is a map $M_1:\Omega\times\mathcal A_1\rightarrow[0,1]$ such that: 
	a) $\forall \omega\in\Omega$, $M_1(\omega,\cdot)$ is a  probability
	measure on
	$\mathcal A_1$; b) $\forall A_1\in\mathcal A_1$, $M_1(\cdot,A_1)$ is an $\mathcal A$-measurable map.
	
  (ii) (Diagonal product of Markov kernels) The diagonal product
	$$M_1\times M_2:(\Omega,\mathcal A)\pt
	(\Omega_1\times\Omega_2,\mathcal A_1\times\mathcal A_2)$$ of two
	Markov kernels $M_1:(\Omega,\mathcal A)\pt   (\Omega_1,\mathcal A_1)$
	and $M_2:(\Omega,\mathcal A)\pt   (\Omega_2,\mathcal A_2)$ is defined
	as the only Markov kernel such that
	$$(M_1\times  M_2)(\omega,A_1\times A_2)=M_1(\omega,A_1)\cdot M_2(\omega,A_2),\quad A_i\in\mathcal A_i,
	i=1,2.
	$$
	
 (iii)
	(Image of a Markov kernel) The image (let us also call it {\it
		probability distribution}) of a Markov kernel $M_1:(\Omega,\mathcal
	A,P)\pt   (\Omega_1,\mathcal A_1)$ on a probability space is the
	probability measure  $P^{M_1}$ on $\mathcal A_1$ defined by
	$P^{M_1}(A_1):=\int_{\Omega}M_1(\omega,A_1)\,dP(\omega)$.
	
(iv) (Independence of Markov kernels, Nogales (2013a))
	Let $(\Omega,\mathcal A,P)$ be a probability space. Two Markov
	kernels $M_1:(\Omega,\mathcal A,P)\pt   (\Omega_1,\mathcal A_1)$ and
	$M_2:(\Omega,\mathcal A,P)\pt   (\Omega_2,\mathcal A_2)$ are said to
	be independent if $P^{M_1\times M_2}=P^{M_1}\times P^{M_2}$. We write $M_1\cind M_2$  (or $M_1\cind_P M_2$).

(v) (Expectation of a Markov kernel) A Markov kernel $M_1:(\Omega,\mathcal A,P)\pt \mathbb R^k$ is said to be $P$-integrable if the map $\omega\mapsto \int_{\mathbb R^k}xM_1(\omega,dx)$ is $P$-integrable, i.e., if there exists and is finite the integral
	$$\int_\Omega\int_{\mathbb R^k}xM_1(\omega,dx)dP(\omega)
	$$
or,  equivalently, if the distribution $(P\otimes M_1)^{\pi_2}$ has finite mean, where $\pi_2:\Omega\times\mathbb R^k\rightarrow\mathbb R^k$ denotes the second coordinatewise projection. In this case, we define the expectation of the Markov kernel $M_1$ as
	$$E_P(M_1):=\int_\Omega\int_{\mathbb R^k}xM_1(\omega,dx)dP(\omega)
	$$
	\end{mydef}

\begin{mydef}\rm
	Let $M_1:(\Omega,\mathcal A,P)\pt \mathbb R^k$ be a $P$-integrable Markov kernel. We define a set function $M_1\cdot P$ on $\mathcal A$ by
	$$(M_1\cdot P)(A):=\int_A\int_{\mathbb R^k}xM_1(\omega,dx)dP(\omega).
	$$
\end{mydef}

Note that $M_1\cdot P\ll P$ and $(M_1\cdot P)^{M_2}\ll P^{M_2}$, when $M_2:(\Omega,\mathcal A,P)\pt (\Omega_2,\mathcal A_2)$ is another Markov kernel.

\begin{mydef}\rm \label{def8} (Conditional expectation of a Markov kernel given another, Nogales (2020)) Let $M_1:(\Omega,\mathcal A,P)\pt \mathbb R^k$ be a $P$-integrable Markov kernel and $M_2:(\Omega,\mathcal A,P)\pt (\Omega_2,\mathcal A_2)$ be a Markov kernel. The conditional expectation $E_P(M_1|M_2)$ is defined by:
	$$E_P(M_1|M_2):=\frac{d(M_1\cdot P)^{M_2}}{dP^{M_2}}
	$$
	i.e., $E_P(M_1|M_2)$ is the (equivalence class of) real measurable function(s) on $(\Omega_2,\mathcal A_2)$ such that, for every $A_2\in\mathcal A_2$,
	\begin{gather*}\begin{split}
			\int_\Omega M_2(\omega,A_2)\int_{\mathbb R^k}xM_1(\omega,dx)dP(\omega)&=\int_{A_2}E_P(M_1|M_2)dP^{M_2}\\
			&=\int_\Omega\int_{A_2}E_P(M_1|M_2)(\omega_2)M_2(\omega,d\omega_2)dP(\omega).
		\end{split}\end{gather*}
	\end{mydef}
	
Several examples and useful remarks and results about the concept above defined can be found in Nogales (2013a), Nogales (2013b) and Nogales (2020). 

\section{Main result}

Let $(\Omega,\mathcal A, P)$ be a probability space, $M:(\Omega,\mathcal A, P)\pt\mathbb R^n$ a Markov kernel with finite mean (i.e. $\int_\Omega\int_{\mathbb R^n}\|x\|_\infty M(\omega,dx)dP(\omega)<\infty$), and $M_i:(\Omega,\mathcal A, P)\pt(\Omega_i,\mathcal A_i)$, $i=1,2$, two arbitrary Markov kernels. 

A previous result will be useful.

\begin{mylemma}\label{Lemma1} \rm $M_2\cind M\times M_1$ if and only if for every bounded functions $f_{01}:(\mathbb R^n\times \Omega_1,\mathcal R^n\times\mathcal A_1)\rightarrow\mathbb R$ and $f_2:(\Omega_2,\mathcal A_2)\rightarrow\mathbb R$ we have that
\begin{gather*}
\int_\Omega\int_{\mathbb R^n\times\Omega_1\times\Omega_2}f_{01}(x,\omega_1)f_2(\omega_2)\,
M\times M_1\times M_2(\omega,d(x,\omega_1,\omega_2))dP(\omega)=\\
\int_\Omega\int_{\mathbb R^n\times\Omega_1}f_{01}(x,\omega_1)\,
M\times M_1(\omega,d(x,\omega_1))dP(\omega) \cdot
\int_\Omega\int_{\Omega_2} f_2(\omega_2)\,
M_2(\omega,d\omega_2)dP(\omega)
\end{gather*}
\end{mylemma}

\begin{myprem} \rm In the statement of the previous lemma we can change {\em bounded\/} by {\em integrable}. 
\end{myprem}

We are now ready for the main result.

\begin{mytheo}\label{theo1} \rm If $M_2\cind M\times M_1$, then 
$E(M|M_1\times M_2)=E(M|M_1)$.
\end{mytheo}

\section{The main theorem in terms of densities}

Let $(\Omega,\mathcal A, P)$ be a probability space and, for $i=1,2,3$, $(\Omega_i,\mathcal A_i,\mu_i)$, $1\le i \le 3$ a $\sigma$-finite measure space and $X_i:(\Omega,\mathcal A, P)\rightarrow (\Omega_i,\mathcal A_i,\mu_i)$ a random variable. Let us consider a fourth random variable $X:(\Omega,\mathcal A, P)\rightarrow (\mathbb R^n,\mathcal R^n,\mu)$ where $\mu$ is the Lebesgue measure or the counter measure on a suitable countable subset of $\mathbb R^n$ when $X$ takes values on it. 

Let us suppose the existence of the next densities: $f_i$ is the $\mu_i$-density of $X_i$, $f_{ij}$ is the $(\mu_i\times\mu_j)$-density of $(X_i,X_j)$, $g_i$ is the $(\mu\times\mu_i)$-density of $(X,X_i)$.

Let us also consider the following Markov kernels:
\begin{gather*}
M=P^{X|X_3}:(\Omega_3,\mathcal A_3,P^{X_3})\pt (\mathbb R^n,\mathcal R^n),\\
M_1=P^{X_1|X_3}:(\Omega_3,\mathcal A_3,P^{X_3})\pt (\Omega_1,\mathcal A_1),\\
M_2=P^{X_2|X_3}:(\Omega_3,\mathcal A_3,P^{X_3})\pt (\Omega_2,\mathcal A_2).
\end{gather*}

It is well known that, for $i=1,2$ and for $P^{X_3}$-almost every $\omega_3$,
$$\phi_i(\omega_3,\omega_i):=\frac{f_{i3}(\omega_i,\omega_3)}{f_3(\omega_3)}\qquad \left(\text{resp.,}\quad \phi(\omega_3,x):=\frac{g_{3}(x,\omega_3)}{f_3(\omega_3)}\right)
$$
is a $\mu_i$-density of $M_i(\omega_3,\cdot)$ (resp. a $\mu$-density of $M(\omega_3,\cdot)$). 

It can be readily shown that, $P^{X_3}$-almost surely, a $(\mu\times\mu_1)$-density of $(M\times M_1)(\omega_3,\cdot)$ is the map
$$(x,\omega_1)\mapsto \phi(\omega_3,x)\cdot\phi_1(\omega_3,\omega_1)=\frac{g_{3}(x,\omega_3)\cdot f_{13}(\omega_1,\omega_3)}{f_3(\omega_3)^2}.
$$

It is shown in Nogales (2013a) that $M_2\cind M\times M_1$ is equivalent to
\begin{gather*}
\int_{\Omega_3}\phi_2(\omega_3,\omega_2)\cdot\phi(\omega_3,x)
\cdot\phi_1(\omega_3,\omega_1)dP^{X_3}(\omega_3)=\\
\int_{\Omega_3}\phi_2(\omega_3,\omega_2)d\mu_3(\omega_3)\cdot
\int_{\Omega_3}\phi(\omega_3,x)
\cdot\phi_1(\omega_3,\omega_1)dP^{X_3}(\omega_3),\quad \mu_2\times\mu\times\mu_1-\text{a.s.}
\end{gather*}
or, which is the same,
\begin{gather*}
	\int_{\Omega_3}\frac{f_{23}(\omega_2,\omega_3)\cdot g_{3}(x,\omega_3)\cdot f_{13}(\omega_1,\omega_3)}{f_3(\omega_3)^2}d\mu_3(\omega_3)=\\
	\int_{\Omega_3}f_{23}(\omega_2,\omega_3)d\mu_3(\omega_3)\cdot
		\int_{\Omega_3}\frac{g_{3}(x,\omega_3)\cdot f_{13}(\omega_1,\omega_3)}{f_3(\omega_3)}d\mu_3(\omega_3),\quad \mu_2\times\mu\times\mu_1-\text{a.s.}
\end{gather*}

In Nogales (2020) it is described how conditional expectations for Markov kernels can be computed when densities are available. In particular,
$$E(M|M_1)(\omega_1)=\int_{\mathbb R^n}x\int_{\Omega_3}\frac{g_{3}(x,\omega_3)\cdot f_{13}(\omega_1,\omega_3)}{f_3(\omega_3)\cdot f_1(\omega_1)}d\mu_3(\omega_3)d\mu(x),\quad P^{X_1}-\text{a.s.}
$$
and
$$E(M|M_1\times M_2)(\omega_1,\omega_2)=\int_{\mathbb R^n}x\int_{\Omega_3}\frac{g_{3}(x,\omega_3)\cdot f_{13}(\omega_1,\omega_3)\cdot f_{23}(\omega_2,\omega_3)}{f_3(\omega_3)\cdot f_1(\omega_1)\cdot f_2(\omega_2)}d\mu_3(\omega_3)d\mu(x), \quad P^{(X_1,X_2)}-\text{a.s.}
$$

\section{An example}

\begin{myexa}\rm 
		Let $\Omega$ be a population with $n$ individuals and consider a partition $(A_{ijkl})_{i,j,k,l=0,1}$ of $\Omega$. We write $n_{ijkl}$ for the number of individuals of $A_{ijkl}$. One or more of the indices $i,j,k,l$ can be replaced by a $+$ sign to denote the union of the corresponding sets of the partition: for instance, $A_{+01+}=A_{0010}\cup A_{1010}\cup A_{0011}\cup A_{1011}$. In particular, $\Omega=A_{++++}$. Similar notations should be used for the numbers $n_{ijkl}$ (e.g. $n_{+0+1}=n_{0001}+n_{0011}+n_{1001}+n_{1011}$). Such a situation will be referred to as 
		\begin{center}S$\left(\ \begin{tabular}{|c|c|c|c|c|c|c|c|}\hline
		$n_{0000}$ & $ n_{0010}$ & $ n_{0100}$ & $ n_{0110}$ & $ n_{1000}$ & $ n_{1010}$ & $ n_{1100}$ & $ n_{1110}$\\ \hline 
		$ n_{0001}$ & $ n_{0011}$ & $ n_{0101}$ & $ n_{0111}$ & $ n_{1001}$ & $ n_{1011}$ & $ n_{1101}$ & $ n_{1111}$\\ \hline
		\end{tabular}\ \right)$
		\end{center}

		We introduce four dichotomic random variables $X_1,X_2,X_3,X$ as follows:
					\begin{gather*} 
					X_1(\omega)=i,\mbox{ if }\omega\in A_{i+++},\  i=0,1,\\
			X_2(\omega)=j,\mbox{ if }\omega\in A_{+j++},\   j=0,1,\\
			X_3(\omega)=k,\mbox{ if }\omega\in A_{++k+},\   k=0,1.\\
			X(\omega)=l,\mbox{ if }\omega\in A_{+++l},\   l=0,1.\\
			\end{gather*}
			
 A such scheme could be obtained when we are interested on the  relationship between two diagnostic procedures, represented by the dichotomous variables $X_1$ and $X_2$ ($X_i=1$ or $0$ when the $i^{th}$ diagnostic test is positive or negative, respectively), for a disease represented by the dichotomous variable $X_3$, which takes the values 1 or 0 depending on whether the disease is actually present or absent.  In this case, we have the following equivalence for some known related concepts:
				\begin{gather*}
				p_3=\mbox{prevalence of the disease}\
				 =\frac{n_{++1+}}{n_{++++}},\\ 
				e_1=\mbox{specificity of }X_1=\frac{n_{+00+}}{n_{++0+}},\ \ e_2=\mbox{specificity of }X_2=\frac{n_{0+0+}}{n_{++0+}},\\ s_1=\mbox{sensitivity of }X_1=\frac{n_{+11}}{n_{++1}},\ \ s_2=\mbox{sensitivity of }X_2=\frac{n_{1+1+}}{n_{++1+}}.
				\end{gather*}	
	The random variable $X$ could represent another disease related in some manner to $X_3$. 
	
	The Markov kernel $M_i:=P^{X_i|X_3}:(\Omega_3,\mathcal A_3)\pt \{0,1\}$, $i=1,2$, can be identified with the matrix
	$$M_i=\left(\begin{array}{cc}
	P(X_i=0|X_3=0) & P(X_i=1|X_3=0)\\
	P(X_i=0|X_3=1) & P(X_i=1|X_3=1)\\
	\end{array}\right)
	=
	\left(\begin{array}{cc}
		e_i & 1-e_i\\
		1-s_i & s_i\\
		\end{array}\right)
	$$
	This way, the distribution of $X_i$ coincides with $(1-p_3,p_3)\cdot M_i$. We also write $M=P^{X|X_3}$ and $Q=P^{X_3}$. 
	
	Note that, for $i,j,k,l=0,1$,
	$$M_1(k,\{i\})=\frac{n_{i+++}}{n_{++k+}},\quad
	M_2(k,\{j\})=\frac{n_{+j++}}{n_{++k+}},\quad
	M(k,\{l\})=\frac{n_{+++l}}{n_{++k+}}.
	$$
	Let us write $M_1(k,i)$ instead of $M_1(k,\{i\})$, for simplicity. We also have
	$$M\times M_1(k,(l,i)):=M(k,l)\cdot M_1(k,i)=
	\frac{n_{+++l}\cdot n_{i+++}}{n_{++k+}^2}.
	$$
	
	The theorem states that
	$$M_2\cind_Q M\times M_1\ \ \Rightarrow\ \  E(M|M_1\times M_2)=E(M|M_1).
	$$
	The statement $M_2\cind_Q M\times M_1$ means that $Q^{M_2\times M\times M_1}=Q^{M_2}\times Q^{M\times M_1}$, that is the same as, for every $i,j,l=0,1$,
	$$\sum_{k=0}^1M_2(k,j)M(k,l)M_1(k,i) P(X_3=k)=P(X_2=j)\sum_{k=0}^1M(k,l)M_1(k,i) P(X_3=k),
	$$
	which is equivalent to
	$$\sum_{k=0}^1\frac{n_{+jk+}\cdot n_{++kl}\cdot n_{i+k+}}{n_{++k+}^2}=
	\frac{n_{+j++}}{n_{++++}}\cdot\sum_{k=0}^1\frac{n_{++kl}\cdot n_{i+k+}}{n_{++k+}}.
	$$
	Writing $\Omega_i=\{0,1\}$, $i=1,2,3$,  $E(M|M_1):\Omega_1\rightarrow\mathbb R$ is defined in such a way that 
	$$\int_{\Omega_3}M_1(\omega_3,A_1)\int_{\mathbb R} x\, M(\omega_3,dx)dQ(\omega_3)=\int_{\Omega_3}\int_{A_1}E(M|M_1)(\omega_1)M_1(\omega_3,d\omega_1)dQ(\omega_3),
	$$
	for every $A_1\subset\{0,1\}$. Taking successively $A_1=\{1\},\, \{0\}$, we get
	\begin{gather*}\begin{split}E(M|M_1=1)&=P(X=1|X_3=0)\cdot P(X_3=0|X_1=1)+P(X=1|X_3=1)\cdot P(X_3=1|X_1=1)\\&=
	\frac{n_{++01}}{n_{++0+}}\cdot \frac{n_{1+0+}}{n_{1+++}}+
	\frac{n_{++11}}{n_{++1+}}\cdot \frac{n_{1+1+}}{n_{1+++}}.
	\end{split}\end{gather*}
	By definition, $E(M|M_1\times M_2):\Omega_1\times\Omega_2\rightarrow\mathbb R$ satisfies
	\begin{gather*}\int_{\Omega_3}M_1(\omega_3,A_1)M_2(\omega_3,A_2)\int_{\mathbb R} x\, M(\omega_3,dx)dQ(\omega_3)=\\
	\int_{\Omega_3}\int_{A_1\times A_2}E(M|M_1\times M_2)(\omega_1,\omega_2)M_1(\omega_3,d(\omega_1,\omega_2))dQ(\omega_3),
	\end{gather*}
	for every $A_1,A_2\subset\{0,1\}$. So, for $i,j=0,1$,
	\begin{gather*}E(M|M_1\times M_2)(i,j)=\\
	\scalebox{0.95}{\includegraphics[width=500pt,height=28pt]{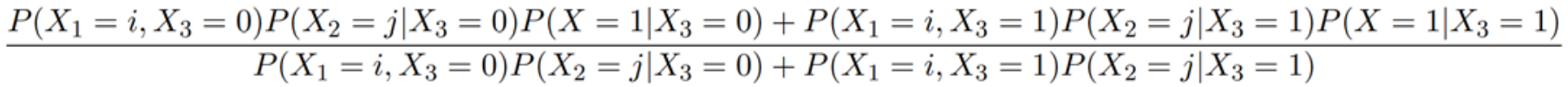}}
	\\=
\frac{n_{i+0+}\cdot  n_{+j0+}\cdot  n_{++1+}\cdot  \frac{n_{++01}}{n_{++0+}} + n_{i+1+}\cdot  n_{+j1+}\cdot  n_{++0+}\cdot  \frac{n_{++11} }{n_{++1+} }}{n_{i+0+}\cdot  n_{+j0+}\cdot  n_{++1+}+ n_{i+1+}\cdot  n_{+j1+}\cdot  n_{++0+}}.
	\end{gather*}

We already have all the necessary ingredients to cook some examples that delimit the main result of the paper. 

	\begin{center}S$\left(\ \begin{tabular}{|c|c|c|c|c|c|c|c|}\hline
	1 & 1 & 2 & 2 & 3 & 3 & 4 & 4\\ \hline 
	3 & 4 & 4 & 3 & 1 & 2 & 2 & 1\\ \hline
	\end{tabular}\ \right)$
\end{center}
is an example where both propositions	$M_2\cind_Q M\times M_1$ and $E(M|M_1\times M_2)=E(M|M_1)$ hold, while		
	\begin{center}S$\left(\ \begin{tabular}{|c|c|c|c|c|c|c|c|}\hline
	1 & 1 & 2 & 2 & 3 & 3 & 4 & 4\\ \hline 
	1 & 2 & 1 & 2 & 3 & 3 & 4 & 4\\ \hline
	\end{tabular}\ \right)$
\end{center}
	is an example where these two sentences fail. Finally, for
	\begin{center}S$\left(\ \begin{tabular}{|c|c|c|c|c|c|c|c|}\hline
	1 & 1 & 2 & 2 & 3 & 3 & 4 & 4\\ \hline 
	3 & 4 & 4 & 1 & 1 & 2 & 2 & 3\\ \hline
	\end{tabular}\ \right),$
\end{center}
    $E(M|M_1\times M_2)=E(M|M_1)$ holds, but not $M_2\cind_Q M\times M_1$.

Note that $E(M)=P(X=1)$ is the prevalence $p$ of the new disease $X$, $E(M_i)=P(X_i=1)$ is the probability of obtaining a positive with the diagnostic procedure $X_i$, and 
$$E(M|M_1)=PPV_1\cdot P(D|D_3)+(1-PPV_1)\cdot P(D|D_3^c)
$$
where $PPV_1$ denotes the positive predictive value for $X_3$ given $X_1$, and $D=\{X=1\}$ and $D_3=\{X_3=1\}$ are the diseased individuals for $X$ and $X_3$, respectively. So, if $X_3$ is considered as a diagnostic procedure for $X$, $P(D|D_3)$ represents the $PPV$ for the disease $X$ given $X_3$ and $P(D|D_3^c)=1-NPV$, where $NPV$ stands for the negative predictive value of $X$ given $X_3$.  $\Box$
		\end{myexa}

\section{Proofs}

{\sc Proof of Lemma 1}.
	The independence $M_2\cind M\times M_1$ is equivalent to $P^{M\times M_1\times M_2}=P^{M\times M_1}\times P^{M_2}$ or, which is the same,
	\begin{gather*}
		\int_\Omega M(\omega,B)\cdot M_1(\omega,A_1)\cdot M_2(\omega,A_2)dP(\omega)=\\
		\int_\Omega M(\omega,B)\cdot M_1(\omega,A_1)dP(\omega)\cdot \int_\Omega M_2(\omega,A_2)dP(\omega)
	\end{gather*}
	for every $B\in\mathcal R^n$ and $A_i\in\mathcal A_i$, $i=1,2$. But this equality can be written as
	\begin{gather*}
		\int_\Omega \int_{\mathbb R^n\times\Omega_1\times\Omega_2}
		I_{B\times A_1}(x,\omega_1)\cdot I_{A_2}(\omega_2)M\times M_1\times M_2(\omega,d(x,\omega_1,\omega_2))dP(\omega)=\\
		\int_\Omega \int_{\mathbb R^n\times\Omega_1}
		I_{B\times A_1}(x,\omega_1)M\times M_1(\omega,d(x,\omega_1))dP(\omega)
		\cdot \int_\Omega \int_{\Omega_2}
		I_{A_2}(\omega_2)M_2(\omega,d\omega_2)dP(\omega),
	\end{gather*}
	and the result follows from here in a standard way. $\Box$
\vspace{2ex}

{\sc Proof of Theorem 1}.
	The conditional expectation $E(M|M_1):(\Omega_1,\mathcal A_1)\rightarrow\mathbb R^n $ is defined in such a way that, for all $A_1\in\mathcal A_1$, 
	\begin{gather*} 
		\int_{\Omega}M_1(\omega,A_1)\int_{\mathbb R^n}xM(\omega,dx)dP(\omega)=
		\int_{A_1}E(M|M_1)dP^{M_1}=\\
		\int_\Omega\int_{A_1}E(M|M_1)(\omega_1)M_1(\omega,d\omega_1)dP(\omega).
	\end{gather*}
	Analogously, $E(M|M_1\times M_2):(\Omega_1\times\Omega_2,\mathcal A_1\times\mathcal A_2)\rightarrow\mathbb R^n$ satisfy
	\begin{gather*} 
		\int_{\Omega}M_1\times M_2(\omega,A_1\times A_2)\int_{\mathbb R^n}xM(\omega,dx)dP(\omega)=
		\int_{A_1\times A_2}E(M|M_1\times M_2)dP^{M_1\times M_2}=\\
		\int_\Omega\int_{A_1\times A_2}E(M|M_1\times M_2)(\omega_1,\omega_2)M_1\times M_2(\omega,d(\omega_1,\omega_2))dP(\omega),
	\end{gather*}
	for every $A_i\in\mathcal A_i$, $i=1,2$. So, it will be enough to prove that, if $M_2\cind M\times M_1$, then
	\begin{gather*} 
		\int_{\Omega}M_1\times M_2(\omega,A_1\times A_2)\int_{\mathbb R^n}xM(\omega,dx)dP(\omega)=\\
		\int_\Omega\int_{A_1\times A_2}E(M|M_1)(\omega_1)\, M_1\times M_2(\omega,d(\omega_1,\omega_2))dP(\omega),
	\end{gather*}
	for every $A_i\in\mathcal A_i$, $i=1,2$. Note that, according to the previous lemma,
	\begin{gather*} 
		\int_\Omega\int_{A_1\times A_2}E(M|M_1)(\omega_1)\, M_1\times M_2(\omega,d(\omega_1,\omega_2))dP(\omega)=\\
		\int_\Omega\int_{\mathbb R^n\times\Omega_1\times\Omega_2}
		I_{A_1}(\omega_1)E(M|M_1)(\omega_1)I_{A_2}(\omega_2)\, M\times M_1\times M_2(\omega,d(x,\omega_1,\omega_2))dP(\omega)=\\
		\int_\Omega\int_{\mathbb R^n\times\Omega_1}
		I_{A_1}(\omega_1)E(M|M_1)(\omega_1)\, M\times M_1(\omega,d(x,\omega_1))dP(\omega)\cdot
		\int_\Omega\int_{\Omega_2}
		I_{A_2}(\omega_2)\, M_2(\omega,d\omega_2)dP(\omega)=\\
		\int_\Omega\int_{A_1}E(M|M_1)(\omega_1)M_1(\omega,d\omega_1)dP(\omega)\cdot \int_\Omega M_2(\omega,A_2)dP(\omega)=\\
		\int_\Omega M_1(\omega,A_1)\int_{\mathbb R^n}xM(\omega,dx)dP(\omega)
		\cdot \int_\Omega M_2(\omega,A_2)dP(\omega)\overset{(*)}{=}\\
		\int_{\Omega}M_1(\omega,A_1)M_2(\omega,A_2)\int_{\mathbb R^n}xM(\omega,dx)dP(\omega),
	\end{gather*}
	where (*) follows from the preceding lemma and the facts that
	$$\int_\Omega M_1(\omega,A_1)\int_{\mathbb R^n}xM(\omega,dx)dP(\omega)=\int_\Omega\int_{\mathbb R^n\times\Omega_1}xI_{A_1}(\omega_1)M\times M_1(\omega,d(x,\omega_1))dP(\omega)
	$$
	and
	$$\int_\Omega M_2(\omega,A_2)dP(\omega)=\int_\Omega\int_{\Omega_2}
	I_{A_2}(\omega_2)\, M_2(\omega,d\omega_2)dP(\omega).\ \ \Box
	$$
	\vspace{2ex}

\section{Acknowledgements}
This paper has been supported by the Junta de Extremadura (Spain) under the grant Gr18016.
\vspace{1ex}

\section* {References:}

\begin{itemize}
	\item Dellacherie, C., Meyer, P.A.:  Probabilities and Potentiel C, North-Holland, Amsterdam (1988).
	
\item Florens, J.P., Mouchart, M., and Rolin, J.M. (1990) Elements of Bayesian Statistics, Marcel Dekker, New York.

	\item Heyer, H.:  Theory of Statistical Experiments, Springer, Berlin (1982).
	
	\item Nogales, A.G.: On Independence of Markov Kernels and a Generalization of Two Theorems of Basu, Journal of Statistical Planning and Inference 143, 603-610  (2013a).
	
	\item Nogales, A.G.: Existence of Regular Conditional Probabilities for Markov kernels, Statistics and Probability Letters 83, 891-897 (2013b).
	
	\item Nogales, A.G.: Conditional Expectation of a Markov Kernel Given Another with some Applications in Statistical Inference and Disease Diagnosis,   Statistics 54 (2), 239--256, (2020).
	
		\item Williams, D.: Probability with Martingales, Cambridge University Press, 1991.
	
\end{itemize}

\end{document}